\documentclass[11pt,a4paper]{article}

\usepackage[UKenglish]{babel}
\usepackage[T1]{fontenc}
\usepackage[utf8]{inputenc}
\usepackage{amsmath,amsthm,amssymb,mathtools}
\usepackage{fullpage}
\usepackage{fourier}
\usepackage{microtype}
\usepackage{enumitem}
\usepackage{needspace}
\usepackage[hidelinks]{hyperref}
\usepackage{pgfplots}
\usepgfplotslibrary{groupplots}
\pgfplotsset{compat=1.18}

\newtheorem{thm}{Theorem}[section]
\newtheorem{lemma}[thm]{Lemma}
\newtheorem{proposition}[thm]{Proposition}
\newtheorem{cor}[thm]{Corollary}
\theoremstyle{remark}
\newtheorem{rem}[thm]{Remark}
\theoremstyle{plain} 
\newtheorem{thmalpha}{Theorem}

\numberwithin{equation}{section}

\newcommand{\R}{\mathbb R}
\newcommand{\N}{\mathbb N}
\newcommand{\Prob}{\mathcal P}
\newcommand{\calE}{\mathcal E}
\newcommand{\PP}{\mathbb P}
\newcommand{\EE}{\mathbb E}
\newcommand{\one}{\mathbf 1}
\newcommand{\dd}{\,\mathrm d}
\newcommand{\weak}{\Rightarrow}
\newcommand{\weakspace}{\Prob([0,\infty))}
\DeclareMathOperator{\vol}{vol}
\DeclareMathOperator{\Law}{Law}

\DeclareMathOperator{\Leb}{Leb}
\DeclareMathOperator{\Exp}{Exp}

\allowdisplaybreaks[2]
\setlist[enumerate]{label=(\roman*),leftmargin=2.2em}
\hypersetup{pdftitle={Large deviations and the asymptotic geometry of Lorentz balls},
  pdfauthor={Simon Matzberger and Joscha Prochno}}

\title{Large deviations and the asymptotic geometry\texorpdfstring{\\}{ }of Lorentz balls}
\author{Simon Matzberger\footnote{University of Passau, Faculty of Computer Science and Mathematics, Dr.-Hans-Kapfinger-Stra{\ss}e 30, 94032 Passau, Germany. Email: matzbe01@ads.uni-passau.de} \and Joscha Prochno\footnote{University of Passau, Faculty of Computer Science and Mathematics, Dr.-Hans-Kapfinger-Stra{\ss}e 30, 94032 Passau, Germany. Email: joscha.prochno@uni-passau.de}}
\date{} 

\begin{document}
\maketitle

\begin{abstract}
We study the uniform distribution on finite-dimensional Lorentz balls
$B_{q,p}^n$ in the range $1<p<q<\infty$. We prove a Sanov-type large
deviation principle for the empirical measures of the absolute coordinates
of a random vector uniformly distributed on $n^{1/q}B_{q,p}^n$.
The rate function is given by an entropy deficit under a weighted quantile
constraint. We show that the entropy maximizer is unique and characterize
its distribution function as the unique increasing solution of a nonlinear
boundary-value problem. The limiting distribution has compact support,
for which we obtain an explicit bound. These results establish the
conjectured empirical limit of Kabluchko, Prochno and Sonnleitner [Journal 
of Functional Analysis, 288(1):110682, 2025]  in the range $1<p<q<\infty$. 
As consequences, we determine the exponential volume asymptotics of Lorentz 
balls and obtain asymptotic independence of any fixed number of coordinates. 
We also prove a law of large numbers for Lorentz norms with finite parameters and 
use it, together with the volume asymptotics, to obtain a
Schechtman--Schmuckenschl\"ager type theorem for intersections of
volume-normalized Lorentz balls. We further derive a variational
formula for the logarithmic intersection asymptotics.
The proof combines an exponential
change of measure with quantile interpolation and a variational argument
that does not assume regularity of the maximizing quantile in advance.
\end{abstract}



\tableofcontents

\section{Introduction and main results}
\label{sec:introduction}

The probabilistic study of high-dimensional convex bodies is a central
part of asymptotic geometric analysis. Endowing a convex body with
normalized Lebesgue measure connects its geometric properties to the
distribution of a random vector and makes probabilistic limit theorems
available for their study. The unit balls of the classical sequence
spaces $\ell_p^n$ provide a particularly instructive setting. A
fundamental tool is the probabilistic representation obtained by
Schechtman and Zinn \cite{SZ1990} and independently by Rachev and
R\"uschendorf \cite{RR1991}: a uniformly distributed point can be
represented using independent $p$-generalized Gaussian random variables,
normalization by their $\ell_p^n$-norm, and an independent radial factor.
This representation underlies the analysis of coordinate distributions,
norms, and asymptotic intersection volumes for these balls.

An important impetus for the systematic study of large deviations in
this setting came from Gantert, Kim and Ramanan \cite{GKR2017}, who
established annealed and quenched large deviation principles for
one-dimensional projections of uniformly distributed points in
$\ell_p^n$-balls. The dependence of the speeds and rate functions on
$p$ reveals geometric information beyond the Gaussian fluctuation
regime. Alonso-Guti\'errez, Prochno and Th\"ale \cite{AGPT2018}
subsequently studied the Euclidean norms of projections onto random
subspaces whose dimensions may grow with the ambient dimension.
Kabluchko, Prochno and Th\"ale \cite{KPT2019} obtained central limit
theorems and large deviation principles for $\ell_q$-norms of random
points in $\ell_p^n$-balls, and later complemented these results by
moderate deviation principles \cite{KPT2021II}, describing the
intermediate scale between central limit and large deviation behaviour.

A closely related development concerns the empirical distribution of
coordinates. Kim and Ramanan \cite{KR2018} proved Sanov-type large
deviation principles for empirical measures under the cone and surface
probability measures on high-dimensional $\ell_p$-spheres and used them
to derive conditional limit theorems. In the matrix setting, Kabluchko,
Prochno and Th\"ale \cite{KPT2020Sanov} established corresponding
large deviation principles for empirical eigenvalue measures of
self-adjoint matrix balls and singular-value measures of Schatten
balls. Such results identify limiting
empirical distributions through variational problems for the rate
function. For an account of these developments and their connections
to geometric functional analysis, we refer to the survey
\cite{Prochno2026Survey}.

The extension of this programme to Orlicz balls brought maximum entropy
methods and Gibbs measures to the fore. Kabluchko and Prochno
\cite{KP2021} related the uniform distribution on an Orlicz ball to
Gibbs distributions with potential given by the defining Orlicz
function. This provided an alternative to the explicit probabilistic
representation for $\ell_p^n$-balls and led to asymptotic volume formulas
(see also \cite{BW2023}) and results on intersection volumes. Building on this approach,
Fr\"uhwirth and Prochno \cite{FP2024} established Sanov-type large
deviation principles for the empirical measures of the coordinates,
together with conditional limit theorems. These results extend the
connection between empirical distributions, entropy maximization, and
asymptotic geometry beyond the classical power constraints.

Lorentz balls form another natural extension of the classical
$\ell_p^n$-balls. For $p<q$, their defining constraint
depends on the order of the absolute coordinates and therefore cannot
be written as a sum of one-coordinate potentials. To make this precise,
let $x_1^*\ge\cdots\ge x_n^*$ denote the decreasing rearrangement of
$(|x_i|)_{i=1}^n$. For $1\le p\le q<\infty$, the finite-dimensional
Lorentz space $\ell_{q,p}^n$ is $\R^n$ equipped with the norm
\[
 \|x\|_{q,p}
 :=\left(\sum_{i=1}^n i^{p/q-1}(x_i^*)^p\right)^{1/p}
 \]
 and we denote its unit ball by
 \[
 B_{q,p}^n:=\{x\in\R^n:\|x\|_{q,p}\le1\}.
\]
When $p=q$, this is the usual $\ell_p^n$-ball. The dependence on the
decreasing rearrangement is the main additional obstacle when $p<q$.

Kabluchko, Prochno and Sonnleitner \cite{KPS2025} developed a
probabilistic approach to Lorentz balls in the case $p=1$. They
obtained an explicit probabilistic representation of the uniform
distribution, identified the limiting empirical law, and derived
coordinate and volumetric limit theorems, including a
Schechtman--Schmuckenschl\"ager type result for intersections with
classical $\ell_r^n$-balls. For general
$1\le p\le q<\infty$, they
formulated a conjecture describing the limiting distribution through
a nonlinear differential equation with boundary conditions. Their
heuristic derivation uses a weighted quantile constraint and a maximum
entropy principle.

In this paper, we establish this conjectured description in the range
$1<p<q<\infty$. We first prove a large deviation principle for the
empirical measures on the positive orthant. We then show that the
associated entropy maximization problem has a unique solution and
derive the differential equation together with its boundary conditions.
A converse argument proves uniqueness of the admissible solution of
the boundary-value problem. In particular, the differential equation
characterizes the limiting law. The same analysis determines the
exponential volume asymptotics and yields a
Poincar\'e--Maxwell--Borel principle for Lorentz balls.
We also prove a law of large numbers for Lorentz norms with
parameters $1\le s\le r<\infty$. Together with the volume
asymptotics, this yields a Schechtman--Schmuckenschl\"ager type
theorem for intersections of two volume-normalized Lorentz balls.
We further determine the exponential decay rate of the
intersection volumes below the threshold.

\subsection{Main results}
Throughout the paper, unless stated otherwise, we fix
$1<p<q<\infty$ and use the normalization
\begin{equation}
 \label{eq:Dn}
 D_n:=n^{1/q}B_{q,p}^n
 =\left\{x\in\R^n:
 \frac1n\sum_{i=1}^n
       \left(\frac in\right)^{p/q-1}(x_i^*)^p\le1\right\}.
\end{equation}
Let $D_n^+:=D_n\cap[0,\infty)^n$, let $U^{(n)}$ be uniformly
distributed on $D_n^+$, and put
\[
 L_n:=\frac1n\sum_{i=1}^n\delta_{U_i^{(n)}}.
\]
We equip $\weakspace$ with the topology of weak convergence. For
$\mu\in\weakspace$, its increasing quantile function is
\[
 Q_\mu(t):=\inf\{x\ge0:\mu([0,x])\ge t\},\qquad 0<t<1.
\]
For $t\in(0,1)$, define
\begin{equation}
 \label{eq:energy}
 w(t):=(1-t)^{p/q-1},\qquad
 \calE(\mu):=\int_0^1 w(t)Q_\mu(t)^p\dd t,\qquad
 K:=\big\{\mu\in\weakspace:\calE(\mu)\le1\big\}.
\end{equation}
Every measure in $K$ has $p$th moment at most one. We write $\Leb$ for Lebesgue measure on $[0,\infty)$ and
$\Leb_A$ for its restriction to a Borel set
$A\subseteq[0,\infty)$. 
For a measure $\mu$ and a measurable map $T$ defined on its
underlying space, we write $T_\#\mu$ for the pushforward of $\mu$
under $T$, defined by
$(T_\#\mu)(A):=\mu(T^{-1}(A))$
for every measurable set $A$ in the target space.
For a probability
measure $\mu$ with finite first moment, we write
\[
 h(\mu) :=
 \begin{cases}
 -\displaystyle\int_0^\infty f(x)\log f(x)\dd x
     &:\, \mu=f\dd x,\\[1mm]
 -\infty&:\,\mu\not\ll\Leb,
 \end{cases}
\]
and put \(0\log0:=0\). This entropy is well defined in $[-\infty,\infty)$. Let us define 
\begin{equation}
 \label{eq:hstar}
 h^*:=\sup_{\mu\in K}h(\mu).
\end{equation}

\begin{thmalpha}[Large deviations and volume asymptotics]
\label{thm:ldp}
Let $1<p<q<\infty$. Then $h^*\in\R$, and the empirical measures
$(L_n)_{n\in\N}$ satisfy a large deviation principle on $\weakspace$
with speed $n$ and good rate function
\begin{equation}
 \label{eq:rate}
 I(\mu):=
 \begin{cases}
 h^*-h(\mu)&:\,\calE(\mu)\le1,\\
 +\infty&:\, \calE(\mu)>1.
 \end{cases}
\end{equation}
Moreover,
\begin{equation}
 \label{eq:volume}
 \lim_{n\to\infty}n^{1/q}\vol_n(B_{q,p}^n)^{1/n}=2e^{h^*}.
\end{equation}
\end{thmalpha}

Thus the exponential volume asymptotics are determined by a
one-dimensional variational problem. Our second theorem identifies
its maximizer and the corresponding empirical limit.

\Needspace{7\baselineskip}
\begin{thmalpha}[The limiting distribution]
\label{thm:limit}
Let $1<p<q<\infty$.
\begin{enumerate}
\item There is a unique measure $\mu_*\in K$ with $h(\mu_*)=h^*$,
and $\calE(\mu_*)=1$. Its quantile $Q=Q_{\mu_*}$ extends to
$Q\in C^1([0,1))$, with $Q(0)=0$, and satisfies
\begin{equation}
 \label{eq:quantile-equation}
 \frac1{Q'(s)}=
 \int_s^1 Q(t)^{p-1}(1-t)^{p/q-1}\dd t,
 \qquad 0\le s<1.
\end{equation}
\item The support of $\mu_*$ is $[0,r_{p,q}]$, where
\begin{equation}
 \label{eq:support-bound}
 0<r_{p,q}\le\left(\frac{p^2}{q-p}\right)^{1/p}.
\end{equation}
Its distribution function $G$ belongs to
$C^1([0,r_{p,q}])\cap C^2((0,r_{p,q}))$ and is the unique increasing
solution, with positive derivative in the interior, of
\begin{equation}
 \label{eq:boundary-problem}
 \begin{gathered}
 G''(x)=-x^{p-1}(1-G(x))^{p/q-1}G'(x),\qquad 0<x<r,\\
 G(0)=0,\qquad \lim_{x\uparrow r}G(x)=1,\qquad
 \lim_{x\uparrow r}G'(x)=0.
 \end{gathered}
\end{equation}
Here uniqueness is among pairs $(r,G)$ with $r\in(0,\infty]$,
$G\in C([0,r))\cap C^2((0,r))$ and $0<G<1$, $G'>0$ on $(0,r)$;
the solution necessarily has $r=r_{p,q}<\infty$. Derivatives at
finite endpoints are understood one-sidedly.
\item The empirical measures $L_n$ converge to $\mu_*$ in probability
with respect to the topology of weak convergence.
If $X^{(n)}$ is uniformly distributed on the full scaled ball $D_n$,
then
\begin{equation}
 \label{eq:signed-limit}
 \frac1n\sum_{i=1}^n\delta_{X_i^{(n)}}
 \xrightarrow{\PP}\nu_{q,p},
 \qquad
 \nu_{q,p}:=\frac12\mu_*+\frac12(-\mathrm{id})_\#\mu_*.
\end{equation}
The measure $\nu_{q,p}$ has the continuous density
\begin{equation}
 \label{eq:symmetric-density}
 f_{q,p}(x)=
 \begin{cases}
 \frac12G'(|x|),&|x|\le r_{p,q},\\
 0,&|x|>r_{p,q}.
 \end{cases}
\end{equation}
\end{enumerate}
\end{thmalpha}

Theorem~\ref{thm:limit} proves the empirical-limit conjecture in
\cite{KPS2025} for $1<p<q<\infty$. The cases $p=1$
and $p=q$ are discussed in that paper and correspond, respectively,
to its explicit Lorentz limit and to the classical $p$-Gaussian law.
The compact support in \eqref{eq:support-bound} is a feature of the
strict inequality $p<q$. Our assertion of uniqueness in
\eqref{eq:boundary-problem} includes the terminal derivative condition;
the differential equation with $G(0)=0$ alone does not determine a
unique solution.

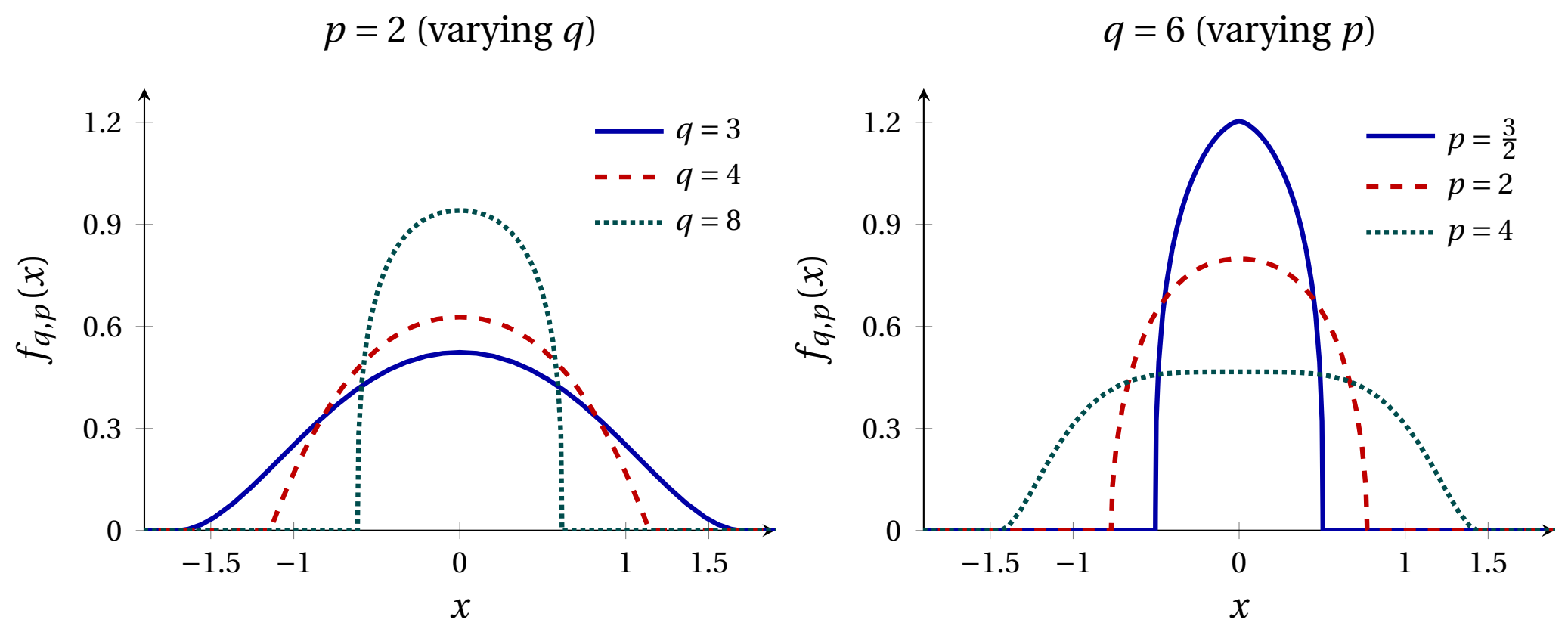
\begin{figure}[t]
\centering
\begin{tikzpicture}
\begin{groupplot}[
    group style={group size=2 by 1,horizontal sep=1.3cm},
    width=7.1cm,
    height=5.45cm,
    xmin=-1.9,xmax=1.9,
    ymin=0,ymax=1.30,
    xtick={-1.5,-1,0,1,1.5},
    ytick={0,0.3,0.6,0.9,1.2},
    tick label style={font=\scriptsize},
    label style={font=\small},
    title style={font=\small},
    xlabel={$x$},
    ylabel={$f_{q,p}(x)$},
    axis lines=left,
    enlargelimits=false,
    legend style={
      font=\scriptsize,
      at={(0.98,0.97)},
      anchor=north east,
      draw=none,
      fill=white,
      fill opacity=0.85,
      text opacity=1,
      row sep=1pt
    },
    legend cell align=left,
]

\nextgroupplot[title={$p=2$ (varying $q$)}]

\addplot+[mark=none, blue!65!black, very thick, solid] coordinates {
(-1.9600000,0.0000000) (-1.7021011,0.0000000) (-1.7019957,0.0000001)
(-1.6999301,0.0000067) (-1.6842265,0.0002809) (-1.6340259,0.0040633)
(-1.5554105,0.0177056) (-1.4786381,0.0385957) (-1.3616882,0.0812517)
(-1.2595616,0.1258570) (-1.1574350,0.1743844) (-1.0553084,0.2244586)
(-0.9531818,0.2740850) (-0.8510551,0.3216207) (-0.7319074,0.3726928)
(-0.6297808,0.4115636) (-0.5276542,0.4451400) (-0.4255276,0.4729428)
(-0.3234010,0.4947028) (-0.2042532,0.5123482) (-0.1021266,0.5209363)
(0.0000000,0.5237282) (0.1021266,0.5209363) (0.2042532,0.5123482)
(0.3234010,0.4947028) (0.4255276,0.4729428) (0.5276542,0.4451400)
(0.6297808,0.4115636) (0.7319074,0.3726928) (0.8510551,0.3216207)
(0.9531818,0.2740850) (1.0553084,0.2244586) (1.1574350,0.1743844)
(1.2595616,0.1258570) (1.3616882,0.0812517) (1.4786381,0.0385957)
(1.5554105,0.0177056) (1.6340259,0.0040633) (1.6842265,0.0002809)
(1.6999301,0.0000067) (1.7019957,0.0000001) (1.7021011,0.0000000)
(1.9600000,0.0000000)
};
\addlegendentry{$q=3$}

\addplot+[mark=none, red!75!black, very thick, dashed] coordinates {
(-1.9600000,0.0000000) (-1.1427671,0.0000080) (-1.1426964,0.0001004)
(-1.1413095,0.0019097) (-1.1307664,0.0155569) (-1.0970623,0.0579217)
(-1.0442811,0.1204815) (-0.9827850,0.1877299) (-0.9142186,0.2558394)
(-0.8570799,0.3072980) (-0.7885136,0.3629774) (-0.7085194,0.4199806)
(-0.6399530,0.4623901) (-0.5713866,0.4991807) (-0.5028202,0.5306641)
(-0.4228261,0.5610952) (-0.3542597,0.5821221) (-0.2856933,0.5987905)
(-0.2171269,0.6113872) (-0.1371328,0.6213107) (-0.0685664,0.6260373)
(0.0000000,0.6275549) (0.0685664,0.6260373) (0.1371328,0.6213107)
(0.2171269,0.6113872) (0.2856933,0.5987905) (0.3542597,0.5821221)
(0.4228261,0.5610952) (0.5028202,0.5306641) (0.5713866,0.4991807)
(0.6399530,0.4623901) (0.7085194,0.4199806) (0.7885136,0.3629774)
(0.8570799,0.3072980) (0.9142186,0.2558394) (0.9827850,0.1877299)
(1.0442811,0.1204815) (1.0970623,0.0579217) (1.1307664,0.0155569)
(1.1413095,0.0019097) (1.1426964,0.0001004) (1.1427671,0.0000080)
(1.9600000,0.0000000)
};
\addlegendentry{$q=4$}

\addplot+[mark=none, teal!60!black, very thick, densely dotted] coordinates {
(-1.9600000,0.0000000) (-0.6155832,0.0002303) (-0.6155451,0.0028755)
(-0.6147980,0.0536458) (-0.6091186,0.2835165) (-0.5909630,0.4343315)
(-0.5625309,0.5534200) (-0.5347654,0.6283336) (-0.4924692,0.7083176)
(-0.4555340,0.7592357) (-0.4185988,0.7989889) (-0.3816636,0.8307321)
(-0.3447284,0.8563725) (-0.3077932,0.8771681) (-0.2647022,0.8964696)
(-0.2277670,0.9094873) (-0.1908318,0.9197547) (-0.1538966,0.9276437)
(-0.1169614,0.9334538) (-0.0738704,0.9379299) (-0.0369352,0.9400241)
(0.0000000,0.9406888) (0.0369352,0.9400241) (0.0738704,0.9379299)
(0.1169614,0.9334538) (0.1538966,0.9276437) (0.1908318,0.9197547)
(0.2277670,0.9094873) (0.2647022,0.8964696) (0.3077932,0.8771681)
(0.3447284,0.8563725) (0.3816636,0.8307321) (0.4185988,0.7989889)
(0.4555340,0.7592357) (0.4924692,0.7083176) (0.5347654,0.6283336)
(0.5625309,0.5534200) (0.5909630,0.4343315) (0.6091186,0.2835165)
(0.6147980,0.0536458) (0.6155451,0.0028755) (0.6155832,0.0002303)
(1.9600000,0.0000000)
};
\addlegendentry{$q=8$}

\nextgroupplot[title={$q=6$ (varying $p$)}]

\addplot+[mark=none, blue!65!black, very thick, solid] coordinates {
(-1.9600000,0.0000000) (-0.5044812,0.0002610) (-0.5044499,0.0032586)
(-0.5038377,0.0607912) (-0.4991834,0.3215984) (-0.4843045,0.4947461)
(-0.4610039,0.6346373) (-0.4382496,0.7252511) (-0.4035871,0.8256630)
(-0.3733181,0.8925857) (-0.3430490,0.9472661) (-0.3127800,0.9931201)
(-0.2825110,1.0321801) (-0.2522419,1.0657588) (-0.2169281,1.0991841)
(-0.1866590,1.1236081) (-0.1563900,1.1445862) (-0.1261210,1.1624129)
(-0.0958519,1.1772708) (-0.0605381,1.1909267) (-0.0302690,1.1993095)
(0.0000000,1.2036723) (0.0302690,1.1993095) (0.0605381,1.1909267)
(0.0958519,1.1772708) (0.1261210,1.1624129) (0.1563900,1.1445862)
(0.1866590,1.1236081) (0.2169281,1.0991841) (0.2522419,1.0657588)
(0.2825110,1.0321801) (0.3127800,0.9931201) (0.3430490,0.9472661)
(0.3733181,0.8925857) (0.4035871,0.8256630) (0.4382496,0.7252511)
(0.4610039,0.6346373) (0.4843045,0.4947461) (0.4991834,0.3215984)
(0.5038377,0.0607912) (0.5044499,0.0032586) (0.5044812,0.0002610)
(1.9600000,0.0000000)
};
\addlegendentry{$p=\tfrac32$}

\addplot+[mark=none, red!75!black, very thick, dashed] coordinates {
(-1.9600000,0.0000000) (-0.7709685,0.0000960) (-0.7709208,0.0011985)
(-0.7699852,0.0223210) (-0.7628722,0.1300775) (-0.7401338,0.2479959)
(-0.7045249,0.3570189) (-0.6697507,0.4322820) (-0.6167781,0.5180381)
(-0.5705198,0.5754723) (-0.5242614,0.6218442) (-0.4780031,0.6598468)
(-0.4317447,0.6911883) (-0.3854863,0.7170439) (-0.3315183,0.7414064)
(-0.2852599,0.7580465) (-0.2390015,0.7712984) (-0.1927432,0.7815638)
(-0.1464848,0.7891756) (-0.0925167,0.7950753) (-0.0462584,0.7978493)
(0.0000000,0.7987327) (0.0462584,0.7978493) (0.0925167,0.7950753)
(0.1464848,0.7891756) (0.1927432,0.7815638) (0.2390015,0.7712984)
(0.2852599,0.7580465) (0.3315183,0.7414064) (0.3854863,0.7170439)
(0.4317447,0.6911883) (0.4780031,0.6598468) (0.5242614,0.6218442)
(0.5705198,0.5754723) (0.6167781,0.5180381) (0.6697507,0.4322820)
(0.7045249,0.3570189) (0.7401338,0.2479959) (0.7628722,0.1300775)
(0.7699852,0.0223210) (0.7709208,0.0011985) (0.7709685,0.0000960)
(1.9600000,0.0000000)
};
\addlegendentry{$p=2$}

\addplot+[mark=none, teal!60!black, very thick, densely dotted] coordinates {
(-1.9600000,0.0000000) (-1.4449653,0.0000000) (-1.4448758,0.0000004)
(-1.4431223,0.0000297) (-1.4297910,0.0010963) (-1.3871741,0.0146712)
(-1.3204352,0.0567353) (-1.2426768,0.1213607) (-1.1559784,0.1967928)
(-1.0837298,0.2548440) (-0.9970314,0.3144518) (-0.8958833,0.3685489)
(-0.8091849,0.4024998) (-0.7224865,0.4267496) (-0.6357881,0.4432491)
(-0.5346400,0.4551943) (-0.4479416,0.4610259) (-0.3612433,0.4641795)
(-0.2745449,0.4656665) (-0.1733968,0.4662711) (-0.0866984,0.4663747)
(0.0000000,0.4663815) (0.0866984,0.4663747) (0.1733968,0.4662711)
(0.2745449,0.4656665) (0.3612433,0.4641795) (0.4479416,0.4610259)
(0.5346400,0.4551943) (0.6357881,0.4432491) (0.7224865,0.4267496)
(0.8091849,0.4024998) (0.8958833,0.3685489) (0.9970314,0.3144518)
(1.0837298,0.2548440) (1.1559784,0.1967928) (1.2426768,0.1213607)
(1.3204352,0.0567353) (1.3871741,0.0146712) (1.4297910,0.0010963)
(1.4431223,0.0000297) (1.4448758,0.0000004) (1.4449653,0.0000000)
(1.9600000,0.0000000)
};
\addlegendentry{$p=4$}

\end{groupplot}
\end{tikzpicture}
\caption{Numerical limiting densities $f_{q,p}$ of the coordinate distributions on
high-dimensional Lorentz balls. The left panel fixes $p=2$ and varies $q$;
the right panel fixes $q=6$ and varies $p$. Each density is symmetric,
has total mass one, and vanishes outside a compact interval.}
\label{fig:lorentz-densities}
\end{figure}

The permutation symmetry of $D_n$ turns convergence of the empirical
measures into convergence of any fixed number of coordinates.

\begin{cor}[A Poincar\'e--Maxwell--Borel principle]
\label{cor:coordinates}
Let $k\in\N$ be fixed and let $X^{(n)}$ be uniformly distributed on
$D_n$. Then
\[
 \Law(X_1^{(n)},\ldots,X_k^{(n)})\weak\nu_{q,p}^{\otimes k}
 \qquad\text{as }n\to\infty.
\]
\end{cor}

We also obtain an intersection-volume result in the sense of
Schechtman and Schmuckenschl\"ager \cite{SS1991}. It complements 
\cite{FP2026, JKP2024, JP2022, KP2021, KPT2020}, extends the intersection 
result of \cite{KPS2025} to $p>1$ and allows a second
Lorentz ball in place of a classical $\ell_r^n$-ball.
For $1\le b\le a<\infty$, write
\begin{equation}
 \label{eq:volume-normalized-lorentz}
 \widehat B_{a,b}^n:=\vol_n(B_{a,b}^n)^{-1/n}B_{a,b}^n,
 \qquad
 c_{a,b}:=\lim_{n\to\infty}
               n^{1/a}\vol_n(B_{a,b}^n)^{1/n}.
\end{equation}
If $1<b<a$, let $h_{a,b}^*$ denote the entropy maximum for
the pair $(a,b)$, so that $h_{q,p}^*=h^*$. Theorem~\ref{thm:ldp},
the volume asymptotics in \cite{KPS2025}, and the classical
$\ell_a^n$-ball volume formula give, respectively,
\begin{equation}
 \label{eq:lorentz-volume-constants}
 c_{a,b}=
 \begin{cases}
  2e^{h_{a,b}^*}&:\,1<b<a,\\
  2e^{1/a}/a&:\,b=1<a,\\
  2\Gamma(1+1/a)(ae)^{1/a}&:\,b=a.
 \end{cases}
\end{equation}

\begin{cor}[Intersection volumes]
\label{cor:intersections}
Let $1<p<q<\infty$ and $1\le s\le r<\infty$, and let
$Q=Q_{\mu_*}$ be the quantile from Theorem~\ref{thm:limit}.
Define
\begin{align}
 m_{q,p;r,s}
 &:=\left(\int_0^1(1-u)^{s/r-1}Q(u)^s\dd u\right)^{1/s},
 \label{eq:cross-lorentz-moment}\\
 \tau_{q,p;r,s}
 &:=\frac{c_{r,s}}{c_{q,p}}m_{q,p;r,s}.
 \label{eq:intersection-threshold}
\end{align}
Then, for every $t>0$ with $t\ne\tau_{q,p;r,s}$,
\begin{equation}
 \label{eq:intersection-limit}
 \lim_{n\to\infty}
 \vol_n\bigl(\widehat B_{q,p}^n\cap t\widehat B_{r,s}^n\bigr)
 =
 \begin{cases}
  0&:\,t<\tau_{q,p;r,s},\\
  1&:\,t>\tau_{q,p;r,s}.
 \end{cases}
\end{equation}
For every fixed $t\in(0,\tau_{q,p;r,s})$, there exists
$J(t)\in(0,\infty)$ such that
\[
 \vol_n\bigl(\widehat B_{q,p}^n\cap t\widehat B_{r,s}^n\bigr)
 =\exp\{-nJ(t)+o(n)\}
 \qquad\text{as }n\to\infty.
\]
Moreover, $\tau_{q,p;r,s}\ge1$. If $1<s<r<\infty$, then
$\tau_{q,p;r,s}=1$ if and only if $(r,s)=(q,p)$.
In particular, the intersection volume at $t=1$ tends to zero
for two distinct parameter pairs in the strict parameter range.
\end{cor}

For $1<s<r$, the threshold can equivalently be written as
$\tau_{q,p;r,s}=e^{h_{r,s}^*-h^*}m_{q,p;r,s}$.
For a classical $\ell_r^n$-ball, corresponding to $s=r$, it is
\begin{equation}
 \label{eq:classical-intersection-threshold}
 \tau_{q,p;r,r}
 =\Gamma(1+1/r)(re)^{1/r}e^{-h^*}
       \left(\int_0^\infty x^r\,\mu_*(\dd x)\right)^{1/r}.
\end{equation} 
The variational characterization of $J(t)$ is given in
\eqref{eq:intersection-rate}, and the exponential asymptotics
are proved in Proposition~\ref{prop:intersection-rate}.

\subsection{Outline of the proofs}

The first step is to express the Lorentz constraint as an integral
against the empirical quantile. The resulting weight depends on $n$,
so the admissible sets vary with the dimension. We control these sets
inside a common weakly compact moment ball and identify all their
possible weak limits. An exponential change of measure then relates
restricted volumes to the empirical measures of independent
exponential random variables. Sanov's theorem gives the upper bound
after truncation of the first-moment tilt. For the lower bound, we use
boundedly supported comparison measures and approximate arbitrary
finite-entropy admissible measures by truncation and dilation.

The variational analysis is carried out in quantile coordinates.
Concavity of the entropy along quantile interpolation and strict
convexity of the energy give uniqueness. Before taking variations, we
remove any initial offset and singular part of the maximizing
quantile. Multiplicative perturbations of its derivative then yield
\eqref{eq:quantile-equation} with an explicit integrable majorant.
The integral equation supplies the required regularity, the boundary
conditions, and the compact-support estimate. Finally, an entropy
comparison shows that every admissible solution of the boundary-value
problem is the same maximizer.

For the intersection results, convexity yields uniform exponential
tails for the one-coordinate marginals. This strengthens empirical
convergence to Wasserstein convergence and gives a law of large
numbers for every Lorentz norm with finite parameters. Combining
this with the volume asymptotics yields the threshold. A second
quantile constraint in the restricted-volume argument gives the
logarithmic intersection asymptotics.

Section~\ref{sec:preliminaries} contains the quantile representation
and the compactness and entropy estimates.
Section~\ref{sec:ldp} proves Theorem~\ref{thm:ldp}.
Section~\ref{sec:maximizer} establishes uniqueness of the entropy
maximizer and the empirical and coordinate limits.
Section~\ref{sec:variation} completes the characterization of the
limiting distribution. Section~\ref{sec:intersections} proves the
intersection-volume results.

\section{Quantiles, compactness and entropy}
\label{sec:preliminaries}

We write $\mu_n\weak\mu$ for weak convergence of probability measures.
The space $\weakspace$ is Polish for this topology; throughout, $d$
denotes any compatible metric. We use the facts that weak convergence
implies $Q_{\mu_n}(t)\to Q_\mu(t)$ at every continuity point of
$Q_\mu$, hence for almost every $t\in(0,1)$, and that
\begin{equation}
 \label{eq:pushforward}
 (Q_\mu)_\#\Leb_{(0,1)}=\mu.
\end{equation}
In particular, integration against $\mu$ can be written as integration
of a function of $Q_\mu$ on $(0,1)$. We refer to
\cite[Appendix~D]{DemboZeitouni1998} for the weak-topology facts used below.

\subsection{The discrete constraint}

For $n\in\N$, we define 
\begin{equation}
 \label{eq:discrete-weight}
 w_n(t)=\left(1-\frac{j-1}{n}\right)^{p/q-1}
 \quad\text{if}\quad\frac{j-1}{n}<t\le\frac jn,
 \qquad
 \calE_n(\mu)=\int_0^1w_n(t)Q_\mu(t)^p\dd t.
\end{equation}
We recall the finite-dimensional quantile identity from
\cite{KPS2025} and include its short proof to fix the mesh
convention.

\begin{lemma}
\label{lem:discrete}
For $x\in[0,\infty)^n$, let $L_x:=n^{-1}\sum_{i=1}^n\delta_{x_i}$.
Then
\begin{equation}
 \label{eq:finite-constraint}
 \calE_n(L_x)=\frac1n\sum_{i=1}^n
       \left(\frac in\right)^{p/q-1}(x_i^*)^p
 \quad
 \text{hence}
 \quad
 x\in D_n^+\ \Longleftrightarrow\ \calE_n(L_x)\le1.
\end{equation}
Moreover, for any \(0<t<1\),
\begin{equation}
 \label{eq:weights}
 1\le w_n(t)\le w(t),\qquad w_n(t)\to w(t),
 \qquad \int_0^1w(t)\dd t=\frac qp.
\end{equation}
\end{lemma}

\begin{proof}
On $((j-1)/n,j/n]$, the increasing quantile of $L_x$ equals
$x_{n-j+1}^*$. This remains true when coordinates coincide.
Integrating over these intervals and replacing $n-j+1$ by $i$
gives \eqref{eq:finite-constraint}. Since $p/q-1\in(-1,0)$,
$w$ is increasing and integrable, and $w_n$ takes its value at the
left endpoint of each interval. This proves \eqref{eq:weights}.
\end{proof}

\subsection{A common compact set}

Let
\[
 C_p:=\left\{\mu\in\weakspace:
                 \int_0^\infty x^p\,\mu(\dd x)\le1\right\},
 \qquad m_1(\mu)=\int_0^\infty x\,\mu(\dd x).
\]
The inequality $w_n,w\ge1$ implies that all the sets
$\{\calE_n\le1\}$ and $K$ are contained in $C_p$.

\begin{lemma}
\label{lem:compactness}
The set $C_p$ is weakly compact and $m_1$ is continuous on $C_p$.
The functionals $\calE_n$ and $\calE$ are weakly lower
semicontinuous. If $\mu_n\weak\mu$, then
\begin{equation}
 \label{eq:moving-liminf}
 \calE(\mu)\le\liminf_{n\to\infty}\calE_n(\mu_n).
\end{equation}
Consequently, $K$ is compact and, for every $\delta>0$,
\begin{equation}
 \label{eq:moving-neighbourhood}
 \{\calE_n\le1\}\subseteq C_p\cap K_\delta
 \quad\text{for all sufficiently large }n,
 \qquad K_\delta=\{\mu:d(\mu,K)\le\delta\}.
\end{equation}
\end{lemma}

\begin{proof}
The $p$th-moment functional is lower semicontinuous, as it is the
supremum of the continuous functionals
$\mu\mapsto\int(x^p\wedge R)\dd\mu$. Thus $C_p$ is closed.
The bound $\mu((R,\infty))\le R^{-p}$ holds uniformly on $C_p$,
so Prohorov's theorem gives compactness. Furthermore,
\begin{equation}
 \label{eq:uniform-tail}
 0\le m_1(\mu)-m_{1,R}(\mu)
 \le\int_{(R,\infty)}x\,\mu(\dd x)\le R^{1-p},
 \qquad \text{where}\quad
 m_{1,R}(\mu)=\int(x\wedge R)\dd\mu.
\end{equation}
Hence the bounded continuous functionals $m_{1,R}$ converge uniformly
to $m_1$ on $C_p$.

If $\mu_n\weak\mu$, quantile convergence and \eqref{eq:weights}
allow us to apply Fatou's lemma to $w_nQ_{\mu_n}^p$, giving
\eqref{eq:moving-liminf}. The same argument with a fixed weight
$w$ or $w_n$ proves lower semicontinuity. In particular, all these
functionals are Borel measurable and $K$ is a closed subset of $C_p$ and hence compact.

If \eqref{eq:moving-neighbourhood} were false, there would be
$n_k\to\infty$ and $\mu_k\in C_p$ with
$\calE_{n_k}(\mu_k)\le1$ and $d(\mu_k,K)>\delta$.
A subsequence converges to some $\mu\in C_p$.
The proof of \eqref{eq:moving-liminf}, applied along this
subsequence, gives $\mu\in K$. Hence
$d(\mu_k,K)\le d(\mu_k,\mu)\to0$,
contradicting $d(\mu_k,K)>\delta$.
\end{proof}

\subsection{Entropy}

Let $\rho(\dd x)=e^{-x}\one_{[0,\infty)}(x)\dd x$ be the
$\Exp(1)$ distribution. We use the relative entropy
\[
 H(\mu\mid\rho)=
 \begin{cases}
 \displaystyle\int\log\!\left(\frac{\dd\mu}{\dd\rho}\right)\dd\mu
     &:\, \mu\ll\rho,\\
 +\infty&:\,\text{otherwise}.
 \end{cases}
\]
It is nonnegative and weakly lower semicontinuous; see
\cite[Section~6.2]{DemboZeitouni1998}.

For every probability measure $\mu$ with finite first moment, the
standard entropy identity gives
\begin{equation}
 \label{eq:entropy-identity}
 h(\mu)=m_1(\mu)-H(\mu\mid\rho),
\end{equation}
with both sides understood in the extended sense. Since $m_1$ is
continuous on $C_p$ by Lemma~\ref{lem:compactness}, the weak lower
semicontinuity of relative entropy implies that $h$ is upper
semicontinuous on $C_p$. Moreover, $h(\mu)\le m_1(\mu)\le1$ for
$\mu\in C_p$. On the other hand, the uniform distribution on $[0,b]$,
with $0<b<(p/q)^{1/p}$, belongs to $K$ and has entropy $\log b$.
Consequently, $h^*\in\R$.

\section{Restricted volumes and large deviations}
\label{sec:ldp}

For a Borel set $A\subseteq\weakspace$, define
\[
 V_n(A):=\vol_n\{x\in D_n^+:L_x\in A\}.
\]
This is well defined because $x\mapsto L_x$ is weakly continuous.
Indeed, $\int f\dd L_x=n^{-1}\sum_i f(x_i)$ is continuous for
every $f\in C_b([0,\infty))$.

\begin{proposition}[Restricted-volume bounds]
\label{prop:volumes}
For every closed $F\subseteq\weakspace$ and every open
$O\subseteq\weakspace$,
\begin{align}
 \limsup_{n\to\infty}\frac1n\log V_n(F)
   &\le\sup_{\mu\in F\cap K}h(\mu),\label{eq:volume-upper}\\
 \liminf_{n\to\infty}\frac1n\log V_n(O)
   &\ge\sup_{\mu\in O\cap K}h(\mu);\label{eq:volume-lower}
\end{align}
we use the conventions $\log0=\sup\varnothing=-\infty$.
\end{proposition}

\begin{proof}
\emph{The upper bound.}
Let $Y_1,Y_2,\ldots$ be independent with law $\rho$ and put
$L_n^Y=n^{-1}\sum_{i=1}^n\delta_{Y_i}$. Since
$\rho^{\otimes n}(\dd x)=e^{-\sum_i x_i}\dd x$ on the positive
orthant, Lemma~\ref{lem:discrete} gives the exact identity
\begin{equation}
 \label{eq:volume-tilt}
 V_n(A)=\EE\!\left[
 e^{nm_1(L_n^Y)}
 \one_{\{L_n^Y\in A,\,\calE_n(L_n^Y)\le1\}}\right].
\end{equation}
By Sanov's theorem \cite[Theorem~6.2.10]{DemboZeitouni1998},
the sequence $(L_n^Y)$ satisfies an LDP for weak convergence,
with speed $n$ and good rate function $H(\cdot\mid\rho)$.

Fix $R,\delta>0$ and set $A_\delta=F\cap C_p\cap K_\delta$.
This set is compact. By \eqref{eq:moving-neighbourhood} and
\eqref{eq:uniform-tail}, for all sufficiently large $n$,
\[
 V_n(F)\le e^{nR^{1-p}}
 \EE\!\left[e^{nm_{1,R}(L_n^Y)}
                  \one_{\{L_n^Y\in A_\delta\}}\right].
\]
The restricted Laplace upper bound for a bounded continuous
functional on a closed set
\cite[Appendix~B.2, pp.~321--322]{Ellis2006} therefore yields
\begin{align*}
 \limsup_{n\to\infty}\frac1n\log V_n(F)
 &\le R^{1-p}+
       \sup_{\mu\in A_\delta}
                \bigl(m_{1,R}(\mu)-H(\mu\mid\rho)\bigr)\\
 &\le R^{1-p}+\sup_{\mu\in A_\delta}h(\mu).
\end{align*}
After letting $R\to\infty$, we consider the compact sets
$A_{1/k}$, which decrease to $F\cap K$. If this intersection is nonempty, upper
semicontinuity of $h$ and compactness imply
\[
 \lim_{k\to\infty}\sup_{\mu\in A_{1/k}}h(\mu)
     =\sup_{\mu\in F\cap K}h(\mu).
\]
For completeness, choose a maximizer in each $A_{1/k}$ and
pass to a convergent subsequence in $C_p$; its limit lies in
$F\cap K$, and upper semicontinuity gives the nontrivial
inequality. This argument also applies when the supremum is
$-\infty$. If $F\cap K$ is empty, compactness implies that
$A_{1/k}$ is empty for some $k$. This proves
\eqref{eq:volume-upper}.

\medskip
\emph{The lower bound for bounded support and strict feasibility.}
Let $\mu=f\dd x\in O$ have support in $[0,M]$, finite entropy,
and $\calE(\mu)<1$. For independent random variables $Z_i$
with law $\mu$, their empirical measures $L_n^Z$ converge
weakly to $\mu$ almost surely. This follows, for example, by
applying the strong law to the indicators of intervals with
rational endpoints and using monotonicity of distribution
functions. Quantile convergence and the bound
\[
 0\le w_n(t)Q_{L_n^Z}(t)^p\le M^pw(t)
\]
give $\calE_n(L_n^Z)\to\calE(\mu)<1$ almost surely by dominated
convergence. Moreover, finite entropy implies $f\log f\in L^1$, so
\[
 -\frac1n\sum_{i=1}^n\log f(Z_i)\longrightarrow h(\mu)
 \qquad\text{almost surely}.
\]
Fix $\varepsilon>0$. Let $B_n$ consist of the vectors $z$
for which $f(z_i)>0$ for all $i$ and
\[
 L_z\in O,\qquad \calE_n(L_z)\le1,\qquad
 -\frac1n\sum_{i=1}^n\log f(z_i)\ge h(\mu)-\varepsilon.
\]
Then $\mu^{\otimes n}(B_n)\to1$, and
\[
 \mu^{\otimes n}(B_n)
   =\int_{B_n}\prod_{i=1}^n f(z_i)\dd z
   \le e^{-n(h(\mu)-\varepsilon)}V_n(O).
\]
Taking logarithms and then letting $\varepsilon\downarrow0$
shows that
\begin{equation}
 \label{eq:bounded-lower}
 \liminf_{n\to\infty}\frac1n\log V_n(O)\ge h(\mu).
\end{equation}

\medskip
\emph{Approximation of an arbitrary admissible measure.}
Let $\mu=f\dd x\in O\cap K$ have finite entropy. For
$z_M=\mu([0,M])>0$, let $\mu_M$ be its conditional law on
$[0,M]$, with density $f_M=z_M^{-1}f\one_{[0,M]}$.
Then $\mu_M\weak\mu$, and
\[
 F_{\mu_M}(x)=\min\{F_\mu(x)/z_M,1\}\ge F_\mu(x).
\]
Consequently, $Q_{\mu_M}\le Q_\mu$ and
$\calE(\mu_M)\le\calE(\mu)\le1$. Also,
\begin{equation}
 \label{eq:truncated-entropy}
 h(\mu_M)=-\frac1{z_M}\int_0^M f\log f\dd x+\log z_M
       \longrightarrow h(\mu).
\end{equation}
Here $z_M\uparrow1$, and $f\log f\in L^1$ by finite entropy.

For $c>0$, write $c\mu$ for the pushforward of $\mu$ by
$x\mapsto cx$. Directly from the density and quantile formulas,
\begin{equation}
 \label{eq:scaling}
 \calE(c\mu)=c^p\calE(\mu),\qquad
 h(c\mu)=h(\mu)+\log c.
\end{equation}
For large $M$ and then $c<1$ sufficiently close to one,
the measure $c\mu_M$ belongs to $O$, has bounded support and
energy strictly below one, and has entropy arbitrarily close
to $h(\mu)$. Apply \eqref{eq:bounded-lower} and take the
supremum over the finite-entropy measures in $O\cap K$.
If there are no such measures, the lower bound is vacuous.
This proves \eqref{eq:volume-lower}.
\end{proof}

\begin{proof}[Proof of Theorem~\ref{thm:ldp}]
Apply Proposition~\ref{prop:volumes} to the whole space to obtain
\begin{equation}
 \label{eq:orthant-volume}
 \lim_{n\to\infty}\frac1n\log\vol_n(D_n^+)=h^*.
\end{equation}
Since $U^{(n)}$ is uniform on $D_n^+$,
\[
 \PP(L_n\in A)=\frac{V_n(A)}{\vol_n(D_n^+)}.
\]
Subtracting \eqref{eq:orthant-volume} from the restricted-volume
bounds gives the LDP upper bound for closed sets and lower
bound for open sets, with rate function \eqref{eq:rate}.
The rate function is nonnegative. For every finite $a$,
\[
 \{I\le a\}=\{\mu\in K:h(\mu)\ge h^*-a\}
\]
is a closed subset of the compact set $K$, by upper
semicontinuity of $h$. Thus $I$ is a good rate function.

Finally, orthant symmetry and \eqref{eq:Dn} give
\[
 2^n\vol_n(D_n^+)=\vol_n(D_n)
       =n^{n/q}\vol_n(B_{q,p}^n).
\]
Taking $n$th roots and using \eqref{eq:orthant-volume}
proves \eqref{eq:volume}.
\end{proof}

\section{The entropy maximizer and empirical limits}
\label{sec:maximizer}

Existence of a maximizer follows from compactness of $K$ and upper
semicontinuity of the entropy. To prove uniqueness, we interpolate
between quantiles. This requires an entropy formula that remains
valid when the quantile has a singular part.

\subsection{Entropy in quantile coordinates}

\begin{lemma}[Quantile entropy formula]
\label{lem:quantile-entropy}
Let $Q:(0,1)\to[0,\infty)$ be nondecreasing and let
$\mu=Q_\#\Leb_{(0,1)}$. Then
\[
 \mu\ll\Leb\quad\Longleftrightarrow\quad Q'>0
                         \quad\text{almost everywhere}.
\]
If these conditions hold and $\mu$ has finite first moment, then
\begin{equation}
 \label{eq:quantile-entropy}
 h(\mu)=\int_0^1\log Q'(t)\dd t
\end{equation}
where both sides are allowed to take the value $-\infty$. In particular, if $h(\mu)\in\R$,
then $\log Q'\in L^1(0,1)$.
\end{lemma}

\begin{proof}
Suppose $\mu=f\dd x$, and let $F$ be its distribution function.
Since $\mu$ has no atoms, $Q$ is strictly increasing: an interval
on which $Q$ is constant would create an atom. The pushforward
identity therefore gives $F(Q(t))=t$ for $0<t<1$.
The equality $F'=f$ holds outside a Lebesgue-null set. Its
preimage under $Q$ has measure zero because $Q_\#\Leb_{(0,1)}=\mu\ll\Leb$.
Since $Q$ is differentiable with finite derivative almost
everywhere, the chain rule gives
\begin{equation}
 \label{eq:inverse-derivative}
 f(Q(t))Q'(t)=1
 \qquad\text{for almost every }t\in(0,1).
\end{equation}
Here $f(Q(t))$ is positive and finite almost everywhere: the
sets where $f=0$ or $f=\infty$ have $\mu$-measure zero after
choosing a finite almost everywhere version of the density.
This proves $Q'>0$ almost everywhere.

Conversely, assume $Q'>0$ almost everywhere. Let $N$ be a
Lebesgue-null Borel set, and let $D$ be the set of points in
$Q^{-1}(N)$ at which $Q'$ exists and lies in $(0,\infty)$.
Then $Q^{-1}(N)\setminus D$ is null. The restriction $Q|_D$
is injective and differentiable relative to $D$. The
one-dimensional area formula on a measurable domain
\cite[Theorem~263D(iv)]{Fremlin2010} gives
\[
 \int_D Q'(t)\dd t=\Leb(Q(D))\le\Leb(N)=0.
\]
Thus $D$ is null and $\mu(N)=\Leb(Q^{-1}(N))=0$, proving
absolute continuity of $\mu$. This argument does not assume
absolute continuity of $Q$.

Under these equivalent conditions, \eqref{eq:inverse-derivative}
holds. Hence $\log Q'=-\log f\circ Q$ almost everywhere.
Apply the pushforward identity to the positive and negative
parts separately. If $\mu$ has finite first moment,
\eqref{eq:entropy-identity} ensures that the entropy is well
defined, and the resulting equality is
\eqref{eq:quantile-entropy}.
\end{proof}

\begin{proposition}[Existence and uniqueness]
\label{prop:unique}
There is a unique entropy maximizer $\mu_*$ in $K$, and
$\calE(\mu_*)=1$.
\end{proposition}

\begin{proof}
Lemma~\ref{lem:compactness}, upper semicontinuity of $h$, and
$h^*\in\R$ give existence and finite maximizing entropy. Every maximizer
has energy one: otherwise, a dilation by some $c>1$
would still belong to $K$ and would increase the entropy
by $\log c$, by \eqref{eq:scaling}.

Suppose that $\mu_0\ne\mu_1$ are two maximizers, with
quantiles $Q_0,Q_1$. For $0<\theta<1$, define
\[
 Q_\theta=(1-\theta)Q_0+\theta Q_1,
 \qquad \mu_\theta=(Q_\theta)_\#\Leb_{(0,1)}.
\]
The function $Q_\theta$ is nondecreasing and left-continuous,
so it is the quantile of $\mu_\theta$. Since $Q_0$ and $Q_1$
differ on a set of positive measure, strict convexity of
$u\mapsto u^p$ and positivity of $w$ give
\begin{equation}
 \label{eq:strict-energy}
 \calE(\mu_\theta)
 <(1-\theta)\calE(\mu_0)+\theta\calE(\mu_1)=1.
\end{equation}

By Lemma~\ref{lem:quantile-entropy}, $Q_0',Q_1'>0$
almost everywhere, and hence $Q_\theta'>0$ almost everywhere.
Since $\mu_\theta\in K\subseteq C_p$, it has finite first
moment. Thus Lemma~\ref{lem:quantile-entropy} and concavity
of the logarithm yield
\[
h(\mu_\theta)
=\int_0^1\log Q_\theta'(t)\dd t
\ge(1-\theta)h(\mu_0)+\theta h(\mu_1)=h^*.
\]
In view of \eqref{eq:strict-energy}, dilating $\mu_\theta$
slightly gives an admissible measure with entropy greater
than $h^*$, a contradiction.

\end{proof}

\subsection{Empirical measures and fixed coordinates}

\begin{proof}[Proof of Theorem~\ref{thm:limit}(iii)]
Fix $\varepsilon>0$ and put
$F_\varepsilon:=\{\mu:d(\mu,\mu_*)\ge\varepsilon\}$.
By compactness, upper semicontinuity, and uniqueness of
the maximizer,
\[
 \sup_{\mu\in F_\varepsilon\cap K}h(\mu)<h^*,
\]
with the same conclusion if the intersection is empty.
Theorem~\ref{thm:ldp} gives
\begin{equation}
 \label{eq:concentration}
 \limsup_{n\to\infty}\frac1n
       \log\PP\bigl(d(L_n,\mu_*)\ge\varepsilon\bigr)
 \le\sup_{\mu\in F_\varepsilon\cap K}h(\mu)-h^*<0.
\end{equation}
This proves $L_n\to\mu_*$ in probability.

To pass to the full scaled ball, use its invariance under
coordinatewise sign changes. We may represent
\[
 X^{(n)}\stackrel d=
 (\varepsilon_1U_1^{(n)},\ldots,\varepsilon_nU_n^{(n)}),
\]
where the $\varepsilon_i$ are independent symmetric signs,
independent of $U^{(n)}$. Orthant boundaries have Lebesgue
measure zero, so this representation is exact in law.
For $f\in C_b(\R)$, put
$\bar f(u)=(f(u)+f(-u))/2$ for $u\ge0$. Conditional on $U^{(n)}$,
the centered sign fluctuation
\[
 Z_n=\frac1n\sum_{i=1}^n
       \bigl(f(\varepsilon_iU_i^{(n)})-\bar f(U_i^{(n)})\bigr)
\]
has mean zero and variance at most $\|f\|_\infty^2/n$.
Hence $Z_n\to0$ in probability, whereas
\[
 \frac1n\sum_{i=1}^n\bar f(U_i^{(n)})
       \xrightarrow{\PP}\int\bar f\dd\mu_*
       =\int f\dd\nu_{q,p}.
\]
Finite intersections of neighborhoods specified by
bounded continuous test functions form a neighborhood
basis for weak convergence. Thus convergence of these
integrals in probability gives \eqref{eq:signed-limit}.
The density formula \eqref{eq:symmetric-density} follows
from the regularity established in
Section~\ref{sec:variation}.
\end{proof}

\begin{proof}[Proof of Corollary~\ref{cor:coordinates}]
Let $M_n=n^{-1}\sum_{i=1}^n\delta_{X_i^{(n)}}$ and fix
$f_1,\ldots,f_k\in C_b(\R)$. Theorem~\ref{thm:limit}
and boundedness give
\[
 \EE\prod_{j=1}^k\int f_j\dd M_n
       \longrightarrow\prod_{j=1}^k\int f_j\dd\nu_{q,p}.
\]
For $n\ge k$, expanding the product gives
\[
\EE\prod_{j=1}^k\int f_j\dd M_n
=\frac{1}{n^k}
\sum_{i_1,\ldots,i_k=1}^n
\EE\prod_{j=1}^k f_j(X_{i_j}^{(n)}).
\] Write
$(n)_k=n(n-1)\cdots(n-k+1)$. By exchangeability, every
term with distinct indices has expectation
$\EE\prod_{j=1}^k f_j(X_j^{(n)})$. The total contribution
of repeated indices therefore gives the estimate
\[
 \left|\EE\prod_{j=1}^k f_j(X_j^{(n)})
        -\EE\prod_{j=1}^k\int f_j\dd M_n\right|
 \le 2\left(1-\frac{(n)_k}{n^k}\right)
                \prod_{j=1}^k\|f_j\|_\infty
 \longrightarrow0.
\]
The one-coordinate marginals converge to $\nu_{q,p}$ and
are tight, hence the $k$-coordinate laws are tight as
well. Every subsequential limit has the displayed
product integrals of $\nu_{q,p}^{\otimes k}$, which
determine this product measure. This proves the claim.
\end{proof}

\section{The variational equation and compact support}
\label{sec:variation}

We now identify the entropy maximizer. The first step
justifies working with the derivative of its quantile.

\begin{lemma}
\label{lem:quantile-ac}
Let $Q=Q_{\mu_*}$. Then $Q'>0$ almost everywhere on $(0,1)$,
and for every $t\in(0,1)$,
\begin{equation}
 \label{eq:quantile-ac}
 Q(t)=\int_0^t Q'(s)\dd s.
\end{equation}
In particular, after setting $Q(0)=0$, the function $Q$ is
absolutely continuous on $[0,b]$ for every $b\in(0,1)$.
\end{lemma}

\begin{proof}
By Lemma~\ref{lem:quantile-entropy}, $v=Q'>0$ almost
everywhere. Monotonicity gives
\[
 0\le\int_a^t v(s)\dd s\le Q(t)-Q(a)\le Q(t),
 \qquad 0<a<t<1.
\]
Consequently, $\widetilde Q(t)=\int_0^t v(s)\dd s$
is finite, satisfies $0\le\widetilde Q\le Q$, and has
derivative $v$ almost everywhere. Since $\widetilde Q$ is continuous and strictly increasing,
it is the quantile function of its induced law $\widetilde\mu$.
The inequality $\widetilde Q\le Q$ implies
$\widetilde\mu\in K$, and the quantile entropy formula gives
\[
 h(\widetilde\mu)=\int_0^1\log v(t)\dd t
                  =h(\mu_*)=h^*.
\]
Uniqueness of the maximizer gives
$\widetilde\mu=\mu_*$. Since $\widetilde Q$ and $Q$
are the corresponding quantiles, they agree on $(0,1)$.
This removes any initial offset, jumps, and singular
continuous part of $Q$.
\end{proof}

\subsection{Admissible variations}

For a quantile $P$, write $\calE(P)$ and $h(P)$ for
the energy and entropy of its induced law. If
$0<\calE(P)<\infty$ and $h(P)\in\R$, scaling gives
\[
 h(P)-\frac1p\log\calE(P)\le h^*.
\]
Since $\log u\le u-1$ for $u>0$, the functional
\begin{equation}
 \label{eq:Phi}
 \Phi(P)=h(P)-\frac1p\calE(P)
\end{equation}
therefore satisfies
\begin{equation}
 \label{eq:Phi-max}
 \Phi(P)\le h(P)-\frac1p\log\calE(P)-\frac1p
        \le h^*-\frac1p=\Phi(Q).
\end{equation}
This observation replaces the constrained problem by
an unconstrained one without introducing an unknown
Lagrange multiplier.

\begin{proposition}
\label{prop:integral-equation}
The maximizing quantile belongs to $C^1([0,1))$ and
satisfies \eqref{eq:quantile-equation}. If
$r=\lim_{s\uparrow1}Q(s)\in(0,\infty]$, then
$G=Q^{-1}$ satisfies
\begin{equation}
 \label{eq:ode}
 G''(x)=-x^{p-1}(1-G(x))^{p/q-1}G'(x),
 \qquad 0<x<r.
\end{equation}
\end{proposition}

\begin{proof}
Write $v=Q'$ and use Lemma~\ref{lem:quantile-ac}.
Let $\psi:(0,1)\to\R$ be bounded and measurable.
For $\varepsilon\in\R$ sufficiently close to zero, define
\[
 v_\varepsilon=v(1+\varepsilon\psi),\qquad
 Q_\varepsilon(t)=\int_0^t v_\varepsilon(s)\dd s.
\]
Put $M=\|\psi\|_\infty$ and assume $|\varepsilon|M<1/2$.
Then $v_\varepsilon>0$ almost everywhere and
\[
 (1-|\varepsilon|M)Q\le Q_\varepsilon
                \le(1+|\varepsilon|M)Q.
\]
Thus $Q_\varepsilon$ induces a measure with finite
positive energy, hence with finite first moment. Moreover,
\[
 \log v_\varepsilon=\log v+\log(1+\varepsilon\psi)
                          \in L^1(0,1),
\]
so it has finite entropy by
Lemma~\ref{lem:quantile-entropy}. Both positive and
negative small $\varepsilon$ are admissible in
\eqref{eq:Phi-max}.

Set $R(t)=\int_0^t v(s)\psi(s)\dd s$, so
$Q_\varepsilon=Q+\varepsilon R$ and $|R|\le MQ$.
Differentiation under the integral sign is justified
by the bounds
\[
 \left|\frac{\psi}{1+\varepsilon\psi}\right|\le2M,
 \qquad
 \bigl|wQ_\varepsilon^{p-1}R\bigr|
       \le M(3/2)^{p-1}wQ^p.
\]
The second majorant is integrable because
$\calE(Q)=1$. Since $\Phi(Q_\varepsilon)$ has a
maximum at zero, its derivative there vanishes:
\begin{equation}
 \label{eq:first-variation}
 0=\int_0^1\psi(s)\dd s
   -\int_0^1w(t)Q(t)^{p-1}
             \left(\int_0^t v(s)\psi(s)\dd s\right)\dd t.
\end{equation}
The double integral is absolutely integrable, since
\[
 \int_0^1\int_0^t
       w(t)Q(t)^{p-1}v(s)|\psi(s)|\dd s\dd t
       \le M\int_0^1w(t)Q(t)^p\dd t=M.
\]
Fubini's theorem turns \eqref{eq:first-variation} into
\[
 \int_0^1\psi(s)
       \left[1-v(s)\int_s^1w(t)Q(t)^{p-1}\dd t\right]\dd s=0.
\]
Because $\psi$ is arbitrary,
\begin{equation}
 \label{eq:a-e-equation}
 v(s)A(s)=1\quad\text{almost everywhere},\qquad
 A(s)=\int_s^1w(t)Q(t)^{p-1}\dd t.
\end{equation}
We now show that $Q\in C^1([0,1))$ and that
\eqref{eq:quantile-equation} holds for every $s\in[0,1)$.
H\"older's inequality gives
\begin{equation}
 \label{eq:A-finite}
 A(0)\le
 \left(\int_0^1wQ^p\right)^{(p-1)/p}
 \left(\int_0^1w\right)^{1/p}
       =\left(\frac qp\right)^{1/p}<\infty.
\end{equation}
Hence $A$ is continuous on $[0,1]$. Since $Q(t)>0$
for $t>0$, we have $A(s)>0$ for $s<1$, whereas
$A(1)=0$. By \eqref{eq:quantile-ac} and
\eqref{eq:a-e-equation},
\[
 Q(t)=\int_0^t\frac{1}{A(s)}\dd s.
\]
This proves $Q\in C^1([0,1))$ and
\eqref{eq:quantile-equation}. On $(0,1)$,
\[
 A'(s)=-w(s)Q(s)^{p-1}.
\]
The inverse function theorem now yields
$G'(x)=A(G(x))$, and therefore
\[
 G''(x)=A'(G(x))G'(x)
       =-x^{p-1}(1-G(x))^{p/q-1}G'(x).
\]
This proves \eqref{eq:ode}.
\end{proof}

\subsection{The support and endpoint conditions}

\begin{proposition}
\label{prop:support}
The support of $\mu_*$ is $[0,r_{p,q}]$, with
$r_{p,q}$ satisfying \eqref{eq:support-bound}.
Its distribution function satisfies
\[
 G(0)=0,\qquad G(r_{p,q})=1,\qquad
 G'(0)>0,\qquad G'(r_{p,q})=0.
\]
Moreover, $G\in C^1([0,r_{p,q}])\cap
C^2((0,r_{p,q}))$, and its density is strictly
decreasing on $(0,r_{p,q})$.
\end{proposition}

\begin{proof}
By monotonicity of $Q$, for $0<s<1$,
\[
 \frac1{Q'(s)}
       =\int_s^1w(t)Q(t)^{p-1}\dd t
       \ge\frac qp Q(s)^{p-1}(1-s)^{p/q}.
\]
Since $Q(s)>0$, it follows that
\[
 (Q(s)^p)'=pQ(s)^{p-1}Q'(s)
       \le\frac{p^2}{q}(1-s)^{-p/q}.
\]
Integrating from zero and using $Q(0)=0$ gives
\begin{equation}
 \label{eq:pointwise-support}
 Q(s)^p\le\frac{p^2}{q-p}
             \bigl(1-(1-s)^{1-p/q}\bigr).
\end{equation}
Thus $r_{p,q}=Q(1-)<\infty$ and
\eqref{eq:support-bound} holds.

The continuous strictly increasing quantile maps
$(0,1)$ onto $(0,r_{p,q})$. Its induced law therefore
has support $[0,r_{p,q}]$, and the inverse $G$ extends
continuously with values zero and one at the endpoints.
Since $G'=A\circ G$ in the interior and $A$ is
continuous on $[0,1]$,
\[
 \lim_{x\downarrow0}G'(x)=A(0)>0,
 \qquad
 \lim_{x\uparrow r_{p,q}}G'(x)=A(1)=0.
\]
The one-sided endpoint derivatives have these values
by the mean value theorem. This proves the stated
regularity. Finally, \eqref{eq:ode} has a strictly
negative right-hand side on $(0,r_{p,q})$, so $G'$
is strictly decreasing there.
\end{proof}

\subsection{Uniqueness of the boundary-value solution}

The preceding results give an admissible solution of
\eqref{eq:boundary-problem}. To obtain a characterization,
we must also show that every such solution yields the
entropy maximizer. The following argument includes
competitors whose quantiles have a singular part.

\begin{proposition}[Converse characterization]
\label{prop:converse}
Let $(r,G)$ satisfy the regularity, positivity and
boundary assumptions in Theorem~\ref{thm:limit}(ii),
as well as \eqref{eq:boundary-problem}. Then
$r=r_{p,q}$ and $G$ is the distribution function
of $\mu_*$.
\end{proposition}

\begin{proof}
Let $Q=G^{-1}$ on $(0,1)$ and set
$A(s)=G'(Q(s))=1/Q'(s)$. The differential equation
gives
\[
 A'(s)=-w(s)Q(s)^{p-1}.
\]
Integrating from $s$ to $b<1$ and letting
$b\uparrow1$, the terminal condition on $G'$
implies
\begin{equation}
 \label{eq:converse-integral}
 A(s)=\int_s^1w(t)Q(t)^{p-1}\dd t,
 \qquad Q'(s)A(s)=1.
\end{equation}
No finite-energy assumption is needed for this
identity: the improper integral equals the finite
quantity $A(s)$.

We have $Q(0+)=0$. The same monotonicity estimate
as in Proposition~\ref{prop:support}, integrated
first from $\varepsilon$ to $s$ and then with
$\varepsilon\downarrow0$, yields
\eqref{eq:pointwise-support}. Consequently,
\[
 r=Q(1-)\le\left(\frac{p^2}{q-p}\right)^{1/p}<\infty.
\]
In particular, $A(0)=\int_0^1wQ^{p-1}$ is finite
and positive, and \eqref{eq:converse-integral}
implies $Q\in C^1([0,1))$ with $Q(0)=0$.
Therefore $Q(t)=\int_0^tQ'(s)\dd s$.
Tonelli's theorem now gives the normalization
\begin{align}
 \calE(Q)
 &=\int_0^1w(t)Q(t)^{p-1}
                   \left(\int_0^tQ'(s)\dd s\right)\dd t\notag\\
 &=\int_0^1 Q'(s)
                   \left(\int_s^1w(t)Q(t)^{p-1}\dd t\right)\dd s
   =1.\label{eq:converse-energy}
\end{align}

The law induced by $Q$ also has finite entropy.
Indeed, $Q'\ge1/A(0)>0$, $Q'$ is bounded on
$[0,1/2]$, and for $s\ge1/2$,
\[
 A(s)\ge\frac qp Q(1/2)^{p-1}(1-s)^{p/q}.
\]
Thus $Q'(s)\le C(1-s)^{-p/q}$ near one, and
$\log Q'\in L^1(0,1)$. The quantile entropy formula from
Lemma~\ref{lem:quantile-entropy} applies.

Let $P$ be the quantile of any $\mu\in K$ with
finite entropy. By monotonicity,
\[
 \int_0^tP'(s)\dd s\le P(t),\qquad 0<t<1,
\]
whether or not $P$ is absolutely continuous.
Using Tonelli, H\"older, and
\eqref{eq:converse-energy}, we obtain
\begin{align*}
 \int_0^1 P'(s)A(s)\dd s
 &=\int_0^1w(t)Q(t)^{p-1}
                   \left(\int_0^tP'(s)\dd s\right)\dd t\\
 &\le\int_0^1w(t)Q(t)^{p-1}P(t)\dd t\\
 &\le\calE(Q)^{(p-1)/p}\calE(P)^{1/p}\le1.
\end{align*}
Finally, $\log u\le u-1$ gives
\begin{equation}
 \label{eq:converse-entropy}
 h(\mu)-h(Q)
       =\int_0^1\log\frac{P'(s)}{Q'(s)}\dd s
       \le\int_0^1P'(s)A(s)\dd s-1\le0.
\end{equation}
Both logarithms are integrable, and the right-hand
side is finite by the preceding estimate, so every
integral in this comparison is justified.
Measures of entropy $-\infty$ do not affect
maximality. Hence the law induced by $Q$ is an
entropy maximizer in $K$, and
Proposition~\ref{prop:unique} identifies it with
$\mu_*$. Its quantile determines $r$ and $G$
uniquely.
\end{proof}

\begin{proof}[Completion of the proof of Theorem~\ref{thm:limit}]
Proposition~\ref{prop:unique} gives the unique
maximizer and energy normalization.
Proposition~\ref{prop:integral-equation} gives
\eqref{eq:quantile-equation} and the differential
equation. Propositions~\ref{prop:support} and
\ref{prop:converse} give the support estimate,
boundary regularity, and uniqueness of the
boundary-value solution. The empirical convergence
was proved in Section~\ref{sec:maximizer}.
The continuous density in
\eqref{eq:symmetric-density} follows from
$G'\in C([0,r_{p,q}])$ and $G'(r_{p,q})=0$.
\end{proof}

\section{Intersection volumes}
\label{sec:intersections}

Throughout this section, $1<p<q<\infty$ are the fixed source
parameters and $1\le s\le r<\infty$ are the target parameters.
We use the volume normalization and constants from
\eqref{eq:volume-normalized-lorentz}--\eqref{eq:lorentz-volume-constants}.
For $0<u<1$, set
\begin{equation}
 \label{eq:target-energy}
 w_{r,s}(u)=(1-u)^{s/r-1},\qquad
 \calE_{r,s}(\mu)=\int_0^1w_{r,s}(u)Q_\mu(u)^s\dd u.
\end{equation}
The corresponding discrete weight and energy are
\begin{equation}
 \label{eq:target-discrete-energy}
 \begin{split}
 w_{n;r,s}(u)&=\left(1-\frac{j-1}{n}\right)^{s/r-1}
       \quad\text{if }\frac{j-1}{n}<u\le\frac jn,\\
 \calE_{n;r,s}(\mu)&=\int_0^1w_{n;r,s}(u)Q_\mu(u)^s\dd u.
 \end{split}
\end{equation}
As in Lemma~\ref{lem:discrete},
\begin{equation}
 \label{eq:target-quantile-identity}
 \calE_{n;r,s}(L_x)=n^{-s/r}\|x\|_{r,s}^s,
 \qquad x\in[0,\infty)^n.
\end{equation}
Also, $1\le w_{n;r,s}\le w_{r,s}$ and
$w_{n;r,s}\to w_{r,s}$ pointwise on $(0,1)$.

\subsection{Moment bounds and Lorentz norms}

We first strengthen the empirical convergence in
Theorem~\ref{thm:limit}. For $k\ge1$, we use the one-dimensional
quantile representation of the $k$-Wasserstein distance,
\[
 W_k(\mu,\nu)
 =\left(\int_0^1|Q_\mu(t)-Q_\nu(t)|^k\dd t\right)^{1/k},
\]
for probability measures with finite $k$th moments.

\begin{lemma}
\label{lem:uniform-tails-wasserstein}
There is a constant $c_0>0$ such that
\begin{equation}
 \label{eq:uniform-coordinate-tail}
 \PP(U_1^{(n)}\ge t)\le2e^{-c_0t},
 \qquad n\in\N,\quad t\ge0.
\end{equation}
Consequently, for every $1\le k<\infty$,
\begin{equation}
 \label{eq:empirical-wasserstein}
 W_k(L_n,\mu_*)\xrightarrow{\PP}0.
\end{equation}
\end{lemma}

\begin{proof}
Permutation symmetry and the bound
$n^{-1}\sum_{i=1}^n x_i^p\le1$ on $D_n^+$ give
\[
 \EE\bigl[(U_1^{(n)})^p\bigr]\le1,
 \qquad \EE U_1^{(n)}\le1.
\]
Put $S_n(t)=\PP(U_1^{(n)}\ge t)$ and
$C_t=D_n^+\cap\{x_1\ge t\}$. Since $D_n^+$ is convex,
\[
 (1-\theta)C_a+\theta C_b
 \subseteq C_{(1-\theta)a+\theta b},
 \qquad 0<\theta<1.
\]
For $S_n(a),S_n(b)>0$, the Brunn--Minkowski inequality,
followed by the weighted arithmetic--geometric mean inequality,
therefore yields
\[
 S_n((1-\theta)a+\theta b)
 \ge\bigl((1-\theta)S_n(a)^{1/n}
             +\theta S_n(b)^{1/n}\bigr)^n
 \ge S_n(a)^{1-\theta}S_n(b)^\theta.
\]
The resulting log-concavity inequality is automatic if either value is zero.
Thus $S_n$ is log-concave. Since $S_n(0)=1$ and
$S_n(2)\le1/2$ by Markov's inequality, log-concavity gives
$S_n(t)\le e^{-c_0t}$ for $t\ge2$, with
$c_0=(\log2)/2$. This proves \eqref{eq:uniform-coordinate-tail}
for all $t\ge0$. In particular,
\begin{equation}
 \label{eq:uniform-coordinate-moments}
 \sup_n\EE\bigl[(U_1^{(n)})^a\bigr]<\infty
 \qquad\text{for every }a>0.
\end{equation}

Fix $k\ge1$ and let $T_R(x)=\min\{x,R\}$.
Write $L_n^R=(T_R)_\#L_n$ and
$\mu_*^R=(T_R)_\#\mu_*$. The natural coupling and
\eqref{eq:uniform-coordinate-moments} imply
\begin{align*}
 \EE W_k(L_n,L_n^R)^k
 &\le\EE\int (x-R)_+^k\,L_n(\dd x)\\
 &\le\frac1R\EE\bigl[(U_1^{(n)})^{k+1}\bigr]
 \le\frac{C_k}{R}.
\end{align*}
For fixed $R$, Theorem~\ref{thm:limit} gives
$L_n^R\to\mu_*^R$ weakly in probability.
On $[0,R]$, weak convergence implies $W_k$-convergence
by almost everywhere quantile convergence and dominated
convergence. Hence $W_k(L_n^R,\mu_*^R)\to0$ in probability.
Moreover, $W_k(\mu_*^R,\mu_*)\to0$ as $R\to\infty$,
since $\mu_*$ has compact support. Apply the triangle
inequality, use Markov's inequality for the first term,
and then let $R\to\infty$ to obtain
\eqref{eq:empirical-wasserstein}.
\end{proof}

\begin{proposition}[A law of large numbers for Lorentz norms]
\label{prop:lorentz-norm-lln}
Let $1\le s\le r<\infty$ and let $X^{(n)}$ be uniformly
distributed on $D_n$. With $m_{q,p;r,s}$ as in
\eqref{eq:cross-lorentz-moment}, we have $0<m_{q,p;r,s}<\infty$ and
\begin{equation}
 \label{eq:cross-lorentz-lln}
 n^{-1/r}\|X^{(n)}\|_{r,s}
 \xrightarrow{\PP}m_{q,p;r,s}.
\end{equation}
\end{proposition}

\begin{proof}
Define
\[
 F_{n;r,s}(\mu)
 =\left(\int_0^1w_{n;r,s}(t)Q_\mu(t)^s\dd t\right)^{1/s}.
\]
We have $w_{n;r,s}\le w_{r,s}$ and pointwise convergence
on $(0,1)$. Choose $k>r$. Since
$w_{r,s}\in L^{k/(k-s)}(0,1)$, Minkowski's and
H\"older's inequalities give, for measures with finite
$k$th moments,
\begin{align*}
 |F_{n;r,s}(\mu)-F_{n;r,s}(\nu)|
 &\le\left(\int_0^1w_{n;r,s}
                  |Q_\mu-Q_\nu|^s\right)^{1/s}\\
 &\le\|w_{r,s}\|_{L^{k/(k-s)}}^{1/s}W_k(\mu,\nu).
\end{align*}
The constant is independent of $n$. Dominated convergence
gives $F_{n;r,s}(\mu_*)\to m_{q,p;r,s}$, because $Q$
is bounded and $w_{r,s}$ is integrable. Thus
Lemma~\ref{lem:uniform-tails-wasserstein} yields
$F_{n;r,s}(L_n)\to m_{q,p;r,s}$ in probability.
By \eqref{eq:target-quantile-identity},
$F_{n;r,s}(L_n)=n^{-1/r}\|U^{(n)}\|_{r,s}$.
The absolute coordinates of $X^{(n)}$ have the law of
$U^{(n)}$, proving \eqref{eq:cross-lorentz-lln}.
Positivity and finiteness of $m_{q,p;r,s}$ follow from
$Q(t)>0$ for $t>0$, boundedness of $Q$, and integrability
of $w_{r,s}$.
\end{proof}

\subsection{The intersection threshold}

\begin{proof}[Proof of Corollary~\ref{cor:intersections}]
Put
\[
 a_n=n^{1/q}\vol_n(B_{q,p}^n)^{1/n},\qquad
 d_n=n^{1/r}\vol_n(B_{r,s}^n)^{1/n}.
\]
Then $a_n\to c_{q,p}$ and $d_n\to c_{r,s}$.
If $X^{(n)}$ is uniformly distributed on $D_n$, then
$X^{(n)}/a_n$ is uniform on $\widehat B_{q,p}^n$.
Since the latter body has volume one,
\begin{equation}
 \label{eq:intersection-probability}
 \vol_n\bigl(\widehat B_{q,p}^n\cap t\widehat B_{r,s}^n\bigr)
 =\PP\!\left(n^{-1/r}\|X^{(n)}\|_{r,s}
                  \le t\frac{a_n}{d_n}\right).
\end{equation}
Proposition~\ref{prop:lorentz-norm-lln} now proves
\eqref{eq:intersection-limit} with the threshold
\eqref{eq:intersection-threshold}.
The bound $\tau_{q,p;r,s}\ge1$ follows from
\[
 \vol_n\bigl(\widehat B_{q,p}^n\cap t\widehat B_{r,s}^n\bigr)
 \le\vol_n(t\widehat B_{r,s}^n)=t^n:
\]
if the threshold were smaller than one, choosing $t$ strictly
between it and one would contradict \eqref{eq:intersection-limit}.

Suppose now that $1<s<r<\infty$, let $\mu_{r,s}$ be the
entropy maximizer for the target pair, and write $m=m_{q,p;r,s}$.
The law of $Y/m$, where $Y\sim\mu_*$, has target energy one
and entropy $h^*-\log m$. Hence
\[
 h_{r,s}^*\ge h^*-\log m.
\]
If $\tau_{q,p;r,s}=1$, equality holds, and uniqueness of the
target maximizer implies
$\mu_{r,s}=(m^{-1}\mathrm{id})_\#\mu_*$.
Thus its distribution function satisfies $G_{r,s}(x)=G(mx)$.
Comparing the differential equations in
Theorem~\ref{thm:limit} for the two parameter pairs and
cancelling the positive derivative $G_{r,s}'(x)$ gives
\[
 m^p x^{p-s}(1-G_{r,s}(x))^{p/q-s/r}=1
\]
on the interior of the target support. Letting $x\downarrow0$
first forces $p=s$ and then $m=1$; the same identity then
gives $q=r$. Conversely, if $(r,s)=(q,p)$, the equality
$\calE(\mu_*)=1$ yields $m_{q,p;q,p}=1$ and
$\tau_{q,p;q,p}=1$. The assertion for $t=1$ follows.
The exponential asymptotics follow from
Proposition~\ref{prop:intersection-rate} below.
\end{proof}

\begin{rem}
\label{rem:critical-intersection}
Corollary~\ref{cor:intersections} does not specify the limit at
$t=\tau_{q,p;r,s}$. A nontrivial critical limit requires control
of the norm fluctuations and of the finer volume asymptotics.
For identical parameter pairs, the threshold equals one and
\[
 \vol_n\bigl(\widehat B_{q,p}^n\cap t\widehat B_{q,p}^n\bigr)
 =\min\{1,t^n\},
\]
so the critical intersection volume equals one in every dimension.
\end{rem}

\subsection{Logarithmic intersection volumes}
 
For $t>0$, we define the following three quantities
\begin{equation}
 \label{eq:intersection-rate}
 \begin{split}
 b_t&=t\frac{c_{q,p}}{c_{r,s}},\\
 K(t)&=\{\mu\in K:\calE_{r,s}(\mu)\le b_t^s\},\\
 J(t)&=h^*-\sup_{\mu\in K(t)}h(\mu).
 \end{split}
\end{equation}
The following proposition identifies the variational quantity
$J(t)$ with the exponential decay rate appearing in
Corollary~\ref{cor:intersections}.

\begin{proposition}[Logarithmic intersection volumes]
\label{prop:intersection-rate}
For every $t>0$, the number $J(t)$ is finite and nonnegative, and
\begin{equation}
 \label{eq:intersection-log-volume}
 \lim_{n\to\infty}\frac1n
 \log\vol_n\bigl(\widehat B_{q,p}^n\cap t\widehat B_{r,s}^n\bigr)
 =-J(t).
\end{equation}
Moreover, $J(t)>0$ for $0<t<\tau_{q,p;r,s}$, whereas
$J(t)=0$ for $t\ge\tau_{q,p;r,s}$.
\end{proposition}

\begin{proof}
Fix $t>0$ and put
\[
 a_n=n^{1/q}\vol_n(B_{q,p}^n)^{1/n},\qquad
 d_n=n^{1/r}\vol_n(B_{r,s}^n)^{1/n},\qquad
 b_n=t\frac{a_n}{d_n}\longrightarrow
 b_t=t\frac{c_{q,p}}{c_{r,s}}.
\]
Define
\[
 \begin{split}
 A_n(t)&=\{\mu:\calE_n(\mu)\le1,
                    \calE_{n;r,s}(\mu)\le b_n^s\},\\
 W_n(t)&=\vol_n\{x\in[0,\infty)^n:L_x\in A_n(t)\}.
 \end{split}
\]
Orthant symmetry and the definitions of $a_n$ and $d_n$ give
\begin{equation}
 \label{eq:ss-volume-ratio}
 \vol_n\bigl(\widehat B_{q,p}^n\cap t\widehat B_{r,s}^n\bigr)
 =\frac{W_n(t)}{\vol_n(D_n^+)}.
\end{equation}
We prove that
\begin{equation}
 \label{eq:ss-positive-volume}
 \lim_{n\to\infty}\frac1n\log W_n(t)
 =\sup_{\mu\in K(t)}h(\mu).
\end{equation}

\emph{The upper bound.}
Since
$1\le w_{n;r,s}\le w_{r,s}$ and $w_{n;r,s}\to w_{r,s}$,
the quantile convergence and Fatou arguments in
Lemma~\ref{lem:compactness} show that $\calE_{r,s}$ and
$\calE_{n;r,s}$ are weakly lower semicontinuous and that
\[
 \calE_{r,s}(\mu)\le
 \liminf_{n\to\infty}\calE_{n;r,s}(\mu_n)
 \qquad\text{whenever }\mu_n\weak\mu.
\]
In particular, all constraints above are Borel measurable and
$K(t)$ is compact. For $\delta>0$, put
\[
 A_\delta(t)=C_p\cap\{\mu:d(\mu,K(t))\le\delta\}.
\]
Then $A_n(t)\subseteq A_\delta(t)$ for all sufficiently large $n$.
Indeed, any sequence contradicting this inclusion has a convergent
subsequence in $C_p$; the two moving lower bounds and $b_n\to b_t$
place its limit in $K(t)$, a contradiction.

Let $Y_i$ be independent with law $\rho$ and write $L_n^Y$ for
their empirical measures. The change of measure in
\eqref{eq:volume-tilt} gives
\[
 W_n(t)=\EE\!\left[e^{nm_1(L_n^Y)}
                         \one_{\{L_n^Y\in A_n(t)\}}\right].
\]
Consequently, for fixed $R,\delta>0$ and all sufficiently large $n$,
\[
 W_n(t)\le e^{nR^{1-p}}
 \EE\!\left[e^{nm_{1,R}(L_n^Y)}
                         \one_{\{L_n^Y\in A_\delta(t)\}}\right].
\]
Sanov's theorem and the restricted Laplace upper bound used in
the proof of Proposition~\ref{prop:volumes} yield
\begin{align*}
 \limsup_{n\to\infty}\frac1n\log W_n(t)
 &\le R^{1-p}+\sup_{\mu\in A_\delta(t)}
                  \bigl(m_{1,R}(\mu)-H(\mu\mid\rho)\bigr)\\
 &\le R^{1-p}+\sup_{\mu\in A_\delta(t)}h(\mu).
\end{align*}
Let $R\to\infty$ and then $\delta\downarrow0$. The compact sets
$A_\delta(t)$ decrease to $K(t)$, so upper semicontinuity of $h$
gives the upper bound in \eqref{eq:ss-positive-volume}.

\emph{The lower bound.}
First let $\nu=f\dd x$ have bounded support, finite entropy, and
\[
 \calE(\nu)<1,\qquad \calE_{r,s}(\nu)<b_t^s.
\]
For independent $Z_i$ with law $\nu$, we have $L_n^Z\weak\nu$
almost surely. Quantile convergence and dominated convergence give
\[
 \calE_n(L_n^Z)\longrightarrow\calE(\nu),\qquad
 \calE_{n;r,s}(L_n^Z)\longrightarrow\calE_{r,s}(\nu)
 \quad\text{almost surely}.
\]
Here both weights $w$ and $w_{r,s}$ are integrable. Thus
$L_n^Z\in A_n(t)$ eventually almost surely. Moreover,
$f\log f\in L^1$, and the strong law gives
$-n^{-1}\sum_{i=1}^n\log f(Z_i)\to h(\nu)$ almost surely.
The density comparison from the lower bound in
Proposition~\ref{prop:volumes} now gives, for every $\varepsilon>0$,
\[
 \begin{split}
 W_n(t)\ge{}& e^{n(h(\nu)-\varepsilon)}
 \PP\!\left(L_n^Z\in A_n(t),
       -\frac1n\sum_{i=1}^n\log f(Z_i)\ge h(\nu)-\varepsilon\right).
 \end{split}
\]
The probability tends to one, hence
$\liminf_n n^{-1}\log W_n(t)\ge h(\nu)$.

For arbitrary $\mu\in K(t)$ with finite entropy, let
$\mu_M=\mu(\,\cdot\mid[0,M])$ whenever $\mu([0,M])>0$.
As in the proof of Proposition~\ref{prop:volumes},
$Q_{\mu_M}\le Q_\mu$ and $h(\mu_M)\to h(\mu)$.
For $0<c<1$, the boundedly supported measure $\nu=c\mu_M$ satisfies
\[
 \calE(\nu)\le c^p<1,\qquad
 \calE_{r,s}(\nu)\le c^sb_t^s<b_t^s,\qquad
 h(\nu)=h(\mu_M)+\log c.
\]
Apply the preceding lower bound, let $M\to\infty$ and $c\uparrow1$,
and take the supremum over $\mu\in K(t)$ of finite entropy.
Such measures exist: the uniform law on a sufficiently short
interval $[0,\eta]$ belongs to $K(t)$ and has entropy $\log\eta$.
Since $K(t)\subseteq C_p$, its entropy supremum is finite.
This proves \eqref{eq:ss-positive-volume}. Combining
\eqref{eq:orthant-volume} and \eqref{eq:ss-volume-ratio} now
proves \eqref{eq:intersection-log-volume}.

The set $K(t)$ is compact, and the entropy is upper
semicontinuous there. Since its supremum is finite and attained,
uniqueness of the entropy maximizer on $K$ shows that
$J(t)=0$ if and only if $\mu_*\in K(t)$. By
\eqref{eq:cross-lorentz-moment} and
\eqref{eq:intersection-threshold}, this is equivalent to
$t\ge\tau_{q,p;r,s}$. This also proves strict positivity
of $J(t)$ below the threshold.
\end{proof}

In particular, if $1<s<r<\infty$ and $(r,s)\ne(q,p)$, then
Corollary~\ref{cor:intersections} gives $\tau_{q,p;r,s}>1$ and
\[
 \vol_n\bigl(\widehat B_{q,p}^n\cap\widehat B_{r,s}^n\bigr)
 =\exp\{-nJ(1)+o(n)\},\qquad J(1)>0.
\]


\section*{Acknowledgement}
The authors were supported by the German Research Foundation under DFG project 516672205 \emph{Limit theorems for the volume of random projections of $\ell_p$-balls}. A considerable part of this work was carried out during a research semester of JP at Graz University of Technology. JP gratefully acknowledges the excellent working conditions and the hospitality provided during his stay. This work is intended to form part of SM's doctoral dissertation.

\section*{Declaration of AI use}
During the preparation of this work, the authors used Claude (Anthropic) and ChatGPT (OpenAI)
as tools for literature search as well as technical and editorial assistance. It was also used to produce Figure \ref{fig:lorentz-densities}. The mathematical framework and proof strategy were fully developed by the authors (and partly rely on \cite{KPS2025}).

\bibliographystyle{abbrv}
\bibliography{references}

\begin{thebibliography}{10}

\bibitem{AGPT2018}
D.~Alonso-Guti{\'e}rrez, J.~Prochno, and C.~Th{\"a}le.
\newblock Large deviations for high-dimensional random projections of
  {$\ell_p^n$}-balls.
\newblock {\em Advances in Applied Mathematics}, 99:1--35, 2018.

\bibitem{BW2023}
F.~Barthe and P.~Wolff.
\newblock Volume properties of high-dimensional {O}rlicz balls.
\newblock In {\em High dimensional probability {IX}---the ethereal volume},
  volume~80 of {\em Progr. Probab.}, pages 75--95. Birkh\"auser/Springer, Cham,
  2023.

\bibitem{DemboZeitouni1998}
A.~Dembo and O.~Zeitouni.
\newblock {\em Large Deviations Techniques and Applications}, volume~38 of {\em
  Applications of Mathematics}.
\newblock Springer, New York, second edition, 1998.

\bibitem{Ellis2006}
R.~S. Ellis.
\newblock {\em Entropy, Large Deviations, and Statistical Mechanics}.
\newblock Classics in Mathematics. Springer, Berlin, 2006.
\newblock Reprint of the 1985 edition.

\bibitem{Fremlin2010}
D.~H. Fremlin.
\newblock {\em Measure Theory. {Vol.} 2: Broad Foundations}.
\newblock Torres Fremlin, Colchester, second edition, 2010.

\bibitem{FP2024}
L.~Fr{\"u}hwirth and J.~Prochno.
\newblock {Sanov}-type large deviations and conditional limit theorems for
  high-dimensional {Orlicz} balls.
\newblock {\em Journal of Mathematical Analysis and Applications},
  536(1):128169, 2024.

\bibitem{FP2026}
L.~Fr\"uhwirth and J.~Prochno.
\newblock Sharp concentration phenomena in high-dimensional {O}rlicz balls.
\newblock {\em J. Funct. Anal.}, 290(6):Paper No. 111322, 45, 2026.

\bibitem{GKR2017}
N.~Gantert, S.~S. Kim, and K.~Ramanan.
\newblock Large deviations for random projections of {$\ell^p$} balls.
\newblock {\em The Annals of Probability}, 45(6B):4419--4476, 2017.

\bibitem{JKP2024}
M.~Juhos, Z.~Kabluchko, and J.~Prochno.
\newblock Limit theorems for mixed-norm sequence spaces with applications to
  volume distribution.
\newblock {\em Electron. J. Probab.}, 29:Paper No. 94, 44, 2024.

\bibitem{JP2022}
M.~Juhos and J.~Prochno.
\newblock Spectral flatness and the volume of intersections of
  {$p$}-ellipsoids.
\newblock {\em J. Complexity}, 70:Paper No. 101617, 21, 2022.

\bibitem{KP2021}
Z.~Kabluchko and J.~Prochno.
\newblock The maximum entropy principle and volumetric properties of {Orlicz}
  balls.
\newblock {\em Journal of Mathematical Analysis and Applications},
  495(1):124687, 2021.

\bibitem{KPS2025}
Z.~Kabluchko, J.~Prochno, and M.~Sonnleitner.
\newblock A probabilistic approach to {Lorentz} balls {$\ell_{q,1}^n$}.
\newblock {\em Journal of Functional Analysis}, 288(1):110682, 2025.

\bibitem{KPT2019}
Z.~Kabluchko, J.~Prochno, and C.~Th{\"a}le.
\newblock High-dimensional limit theorems for random vectors in
  {$\ell_p^n$}-balls.
\newblock {\em Communications in Contemporary Mathematics}, 21(1):1750092,
  2019.

\bibitem{KPT2020}
Z.~Kabluchko, J.~Prochno, and C.~Th\"ale.
\newblock Intersection of unit balls in classical matrix ensembles.
\newblock {\em Israel J. Math.}, 239(1):129--172, 2020.

\bibitem{KPT2020Sanov}
Z.~Kabluchko, J.~Prochno, and C.~Th{\"a}le.
\newblock {Sanov}-type large deviations in {Schatten} classes.
\newblock {\em Annales de l'Institut Henri Poincar{\'e}, Probabilit{\'e}s et
  Statistiques}, 56(2):928--953, 2020.

\bibitem{KPT2021II}
Z.~Kabluchko, J.~Prochno, and C.~Th{\"a}le.
\newblock High-dimensional limit theorems for random vectors in
  {$\ell_p^n$}-balls. {II}.
\newblock {\em Communications in Contemporary Mathematics}, 23(3):1950073,
  2021.

\bibitem{KR2018}
S.~S. Kim and K.~Ramanan.
\newblock A conditional limit theorem for high-dimensional {$\ell^p$}-spheres.
\newblock {\em Journal of Applied Probability}, 55(4):1060--1077, 2018.

\bibitem{Prochno2026Survey}
J.~Prochno.
\newblock The large and moderate deviations approach in geometric functional
  analysis.
\newblock In R.~Adamczak, N.~Gozlan, K.~Lounici, M.~Madiman, F.~Merlev{\`e}de,
  and E.~Werner, editors, {\em High Dimensional Probability X: The B{\k e}dlewo
  Volume}, volume~82 of {\em Progress in Probability}, pages 347--421.
  Birkh{\"a}user, Cham, 2026.

\bibitem{RR1991}
S.~T. Rachev and L.~R{\"u}schendorf.
\newblock Approximate independence of distributions on spheres and their
  stability properties.
\newblock {\em The Annals of Probability}, 19(3):1311--1337, 1991.

\bibitem{SS1991}
G.~Schechtman and M.~Schmuckenschl{\"a}ger.
\newblock Another remark on the volume of the intersection of two {$L_p^n$}
  balls.
\newblock In J.~Lindenstrauss and V.~D. Milman, editors, {\em Geometric Aspects
  of Functional Analysis: Israel Seminar ({GAFA}) 1989--90}, volume 1469 of
  {\em Lecture Notes in Mathematics}, pages 174--178. Springer, Berlin,
  Heidelberg, 1991.

\bibitem{SZ1990}
G.~Schechtman and J.~Zinn.
\newblock On the volume of the intersection of two {$L_p^n$} balls.
\newblock {\em Proceedings of the American Mathematical Society},
  110(1):217--224, 1990.

\end{thebibliography}

\end{document}